\documentclass[12pt]{amsart}
\usepackage{}
\usepackage{amsmath,amssymb,txfonts}
\usepackage{amssymb}
\usepackage{amsxtra}
\usepackage{amsthm, color}
\usepackage{txfonts}
\usepackage{graphicx}

\usepackage{mathrsfs}
\usepackage{graphicx}

\usepackage[active]{srcltx}

\usepackage[italian, english]{babel}
\usepackage{graphicx}
\usepackage[T1]{fontenc}
\usepackage[latin1]{inputenc}
\usepackage{times}
\usepackage{citeref}
\newtheorem{thm}{Theorem}[section]
\newtheorem{lemma}[thm]{Lemma}
\newtheorem{defn}[thm]{Definition}
\newtheorem{corollary}[thm]{Corollary}
\newtheorem{proposition}[thm]{Proposition}
\newtheorem{remark}[thm]{Remark}

    \newcommand{\ct}[1]{< {#1}\rangle \lower.3ex\hbox{$_{t}$}}
    \newcommand{\lt}[1]{[ {#1}] \lower.3ex\hbox{$_{t}$}}

\newcommand{\R}{{\mathbb{R}}}

\newcommand{\sfee}{{\mathbb S}^{n-2}}

\begin{document}

\selectlanguage{english}
\title[Affine Poisson \& Non-Poisson trace principles for $\dot{H}^{-1<-\alpha\le-\frac12}(\R^{n-1})$]{\small Affine Poisson \& Non-Poisson trace principles for $\dot{H}^{-1<-\alpha\le -\frac12}(\R^{n-1\ge 2})$}
\author[N. Lombardi]{N. Lombardi}
\address{Dipartimento di Matematica e Informatica,
Universit\`a di Firenze,	
50134, Firenze, Italy}
\email{
nico.lombardi@unifi.it}
\author[J. Xiao]{J. Xiao}
\address{Department of Mathematics and Statistics,
	Memorial University, St. John's, NL A1C 5S7, Canada}
\email{jxiao@math.mun.ca}
\date{}

\thanks{The project was basically completed in the first half year of 2018 when the first-named author studied at Memorial University as a visiting graduate student supported by the second-named author's NSERC-Discovery Grants Program.}
\maketitle

\begin{abstract}
	This note discovers not only an affine non-sharp-Poisson trace inequality but also its sharp-non-Poisson version for a Sobolev function with the fractional antiderivative.  
\end{abstract}

\section{Introduction}\label{s1}

The Sobolev trace inequalities play a fundamental role in Analysis and Geometry as well as PDE. Especially, the sharp $L^2$ inequality for the half-space $\R^{n\ge 3}_{+}=\R_+ \times \R^{n-1}$ says that for a function $f\colon \R^{n}_{+} \to \R$, being smooth and decaying fast at infinity, one has
\begin{equation}\label{intro1}
\Bigg( \int_{\R^{n-1}}|f(0,x)|^{\frac{2(n-1)}{n-2}} dx  \Bigg)^{\frac{n-2}{n-1}} \leq \mathsf{A}(n) \int_{\R^{n}_{+}} |\nabla f(t,x)|^2 dxdt,
\end{equation}
where
$$\mathsf{A}(n)=\frac{1}{\sqrt{\pi}(n-2)}\Bigg( \frac{\Gamma(n)}{(n-1)\Gamma(\frac{n-1}{2})}  \Bigg)^\frac1{n-1}$$
is the sharp constant and $\nabla$ is the full gradient in $(t,x)$.

This result was proved by Beckner \cite{Beckner} and Escobar \cite{Escobar} independently.
Although both of them employed the conformal invariance of the inequality to link this problem with that one over the unit sphere $\mathbb{S}^{n-2}$ of $\R^{n-1}$, their proofs were quite different. Escobar used the Obata's rigidity theorem, while Beckner deduced his result from the spectral form of Lieb's sharp Hardy-Littlewood-Sobolev inequality.
They also proved that the equality in \eqref{intro1} holds if and only if
$$
f(t,x)=\gamma ((t+\delta)^2 + |x-x_0|^2)^{-\frac{n-2}{2}}
$$
for some
$$(\gamma,\delta,x_0) \in \R\times\R_+\times \R^{n-1}.
$$
We also refer to the Lions work \cite{Lions} on finding the extremal functions of \eqref{intro1} through a concentration-compactness principle.

Moreover, \cite{Xiao} and \cite{Einav/Loss} contain some extended counterparts of \eqref{intro1} for the fractional homogeneous Sobolev space. For $-1<\alpha<1$ let $\dot{H}^{\alpha}(\R^{n-1})$ be the closure of $C_{c}^{\infty}(\R^{n-1})$ (the space of all smooth functions compactly supported in $\R^{n-1}$) under the fractional homogeneous Sobolev norm
$$\|f\|_{\dot{H}^{\alpha}(\R^{n-1})}=
\bigg( \int_{\R^{n-1}} (2\pi |x|)^{2\alpha} |\widehat{f}(x)|^2 dx \bigg)^{\frac{1}{2}}.$$ Note that $$\big\{\dot{H}^{\alpha}(\R^{n-1}), \dot{H}^{-\alpha}(\R^{n-1})\big\}
$$
is a dual pair. So, $ \dot{H}^{-\alpha}(\R^{n-1})$ is called the anti-Sobolev space corresponding to $\dot{H}^{\alpha}(\R^{n-1})$.

In \cite{Xiao} Xiao proved a sharp fractional Sobolev trace inequality in the sense of Beckner and Escobar for a particular class of functions on $\R^n_+$. To be more precise, let $f\colon \R^{n-1} \to \R$ be a function in $\dot{H}^{\alpha}(\R^{n-1})$ and consider its Poisson extension onto $\R_{+}^{n}$, i.e., 
\begin{equation}
\label{P1}
\begin{cases}
f(t,x)=P_t\ast f(x);\\
f(0,x)=\lim_{t \rightarrow 0^+}P_t\ast f(x)=f(x);\\
P_t(x)=t^{1-n}P(t^{-1}x)={\pi^{-\frac{n}{2}}\Gamma(\frac{n}{2})t}{(t^2 + |x|^2)^{-\frac{n}{2}}};\\
P(x)=\pi^{-\frac{n}{2}}\Gamma(\frac{n}{2})(1+|x|^2)^{-\frac{n}{2}},
	\end{cases}
\end{equation}
where $\Gamma$ is the standard gamma function.

Regarding this case Xiao proved in \cite[Theorem 1]{Xiao} that if 
$$
0<\alpha<1\quad\&\quad
f\in\dot{H}^{\alpha}(\R^{n-1}),
$$ then
\begin{equation}\label{intro2}
\Bigg( \int_{\R^{n-1}} |f(x)|^{\frac{2(n-1)}{n-1-2 \alpha}} dx\Bigg)^{\frac{n-1-2 \alpha}{n-1}} \leq \mathsf{B}(n-1, \alpha) \int_{\R^{n}_{+}} |\nabla f(t,x)|^2\frac{dxdt}{t^{2\alpha-1}},
\end{equation}
where 
$$\mathsf{B}(n,\alpha)=\Bigg( \frac{2^{1-4\alpha}}{\pi^{\alpha} \Gamma(2(1-\alpha))}\Bigg) \Bigg( \frac{\Gamma(\frac{n}{2}-\alpha)}{\Gamma(\frac{n}{2}+\alpha)}\Bigg) \Bigg( \frac{\Gamma(n)}{\Gamma(\frac{n}{2})}\Bigg)^{\frac{2\alpha}{n}}
$$
is optimal. We refer also to the work made by Einav and Loss in \cite{Einav/Loss}, where they proved a sharp fractional Sobolev trace inequality in a more general setting. In short - given a function $f\in \dot{H}^{\alpha}(\R^{n})$, instead of taking its Poisson's extension onto $\R_{+}^{n+1}$ they considered its restriction to $\R^{n-m}$, in terms of the trace of $f$, where $m$ is a positive integer obeying 
$0\leq\frac{m}{2}<\alpha<\frac{n}{2}$. 

\vspace{0.25cm} 

The aim of this paper is to obtain a precise dual version of (\ref{intro2}). Still working on the Poisson extension of an arbitrary function 
$f\in \dot{H}^{\alpha}(\R^{n-1})$ we fortunately find a possibility to shrink the right-hand-integral of (\ref{intro2}):
$$
\int_{\R^{n}_{+}} |\nabla f(t,x)|^2 \frac{dxdt}{t^{2\alpha-1}} =\int_{\R^n_+}\Bigg(\bigg(\frac{\partial f}{\partial t}(t,x) \bigg)^2 + \sum_{i=1}^{n-1} \bigg( \frac{\partial f}{\partial x_i}(t,x) \bigg)^2\Bigg)\,{t^{1-2\alpha}}\,dxdt
$$
in accordance with the following elementary rule: if $\nabla_x f(t,x)$  stands for the gradient of $f$ at the point $(t,x)$ with respect to the variable $x\in \R^{n-1}$ only, then
\begin{align*}&\int_{\R^{n}_{+}} |\nabla f(t,x)|^2 t^{1-2\alpha} dxdt\\
 &\ \ =\int_{\R^{n}_{+}} \bigg|\frac{\partial f}{\partial t}(t,x)\bigg|^2  t^{1-2\alpha} dxdt+\int_{\R^{n}_{+}} |\nabla_x f(t,x)|^2 t^{1-2\alpha} dxdt\\
&\ \ \ge 2\Bigg(\int_{\R^{n}_{+}} \bigg|\frac{\partial f}{\partial t}(t,x)\bigg|^2  t^{1-2\alpha} dxdt\Bigg)^{\frac{1}{2}}\Bigg(\int_{\R^{n}_{+}} |\nabla_x f(t,x)|^2 t^{1-2\alpha} dxdt\Bigg)^{\frac{1}{2}},
\end{align*}
and hence the energy of $f$ is split in two terms - one depending on the derivative in $t$ and the other in $x$.
Actually, the above splitting idea suggests us to discover an affine and fractional counterpart of the Escobar and Beckner result. To justify this affine counterpart, we introduce briefly the work made by De N\'apoli, Haddad, Jim\'enez and Montenegro in \cite{DeNapoli}, where they proved a stronger version of \eqref{intro1}. They established a sharp affine $L^2$ Sobolev trace inequality involving the $L^2$-affine energy
$$\mathcal{E}_{2}(f)=c_{n-1} \Bigg( \int_{\sfee} \|\nabla_{\xi} f(t,x)\|_{L^2(\R^{n}_+)}^{1-n}d\xi \Bigg)^{\frac1{1-n}},
$$ 
where 
$$c_{n}=n^{\frac{n+2}{2n}}\sqrt[n]{\omega_{n}},
$$ 
and $\omega_{n}$ denotes the volume of the Euclidean ball $\mathbb{B}^n$ in $\R^n$:
$$\omega_{n}=\frac{\pi^{\frac{n}{2}}}{\Gamma(\frac{n}{2}+1)},$$
and it can be extended to every positive real number $s>0$ via
$$\omega_{s}=\frac{\pi^\frac{s}{2}}{\Gamma(\frac{s+2}{2})}.$$
Moreover we have that for $\xi \in \sfee$ - the unit sphere of $\R^{n-1}$, and $\nabla_{\xi} f(t,x)$ is the directional derivative with respect to the direction $\xi$ at the point $(t,x)$, and 
$$\nabla_{\xi} f(t,x)=\langle\nabla_x f(t,x), \xi\rangle.$$

\noindent
Remarkably, the $L^2$-affine energy and its extension to $p\in (1,n)$ have appeared in many results concerning the affine Sobolev inequalities (cf. \cite{Lutwak, Zhang}) and the affine version of the P\'olya-Szeg\"o principle (cf. \cite{Cianchi, Haberl}). We also remark that the $L^2$ (and $L^p$)-affine energy is invariant under volume preserving affine transformations of $\R_{+}^n$.

In \cite[Theorem 1]{DeNapoli}, De N\'apoli, Haddad, Jim\'enez and Montenegro recently showed that for any $f\in C^\infty_c(\mathbb R\times\R^{n-1})$ one has
\begin{equation}\label{intro3}
\Bigg( \int_{\R^{n-1}} |f(0,x)|^{\frac{2(n-1)}{n-2}} dx  \Bigg)^{\frac{n-2}{n-1}}\leq 2 \mathsf{A}(n) \mathcal{E}_{2}(f)
\bigg\| \frac{\partial f}{\partial t} \bigg\|_{L^2(\R_{+}^n)},
\end{equation}
and this is invariant under the affine transformations of $\R_{+}^n$. Moreover, the equality of \eqref{intro3} holds if and only if 
$$f(t,x)=\pm \big((\lambda t+ \delta)^2 +|B(x-x_0)|^2\big)^{-\frac{n-2}{2}}
$$
for some triple
$$(\lambda, \delta,x_0)\in \R_+\times\R_+\times \R^{n-1}$$ 
and a matrix  $B\in GL_{n-1}$ - the set of all invertible real 
$(n-1)\times (n-1)$-matrices.

We note also that they proved in \cite{DeNapoli} an $L^p$ version of sharp affine Sobolev trace inequality for $p\in (1,n)$. Indeed, in his paper \cite{Escobar} Escobar conjectured a possible extension of his result \eqref{intro1} (verified then by Nazareth in \cite{Nazaret}). De N\'apoli, Haddad, Jim\'enez and Montenegro proved an affine and stronger version of results in \cite{Beckner}, \cite{Escobar} and \cite{Nazaret}. With stronger we mean, for instance in the case $p=2$, that 
$$\mathcal{E}_{2}(f)\leq \bigg(\int_{\R^{n}_{+}} |\nabla_{x}f(t,x)|^2 dxdt\bigg)^{\frac{1}{2}},$$
as we can find in \cite{DeNapoli} and \cite{Lutwak}. Consequently, we have
\begin{align*}&\Bigg( \int_{\R^{n-1}} |f(0,x)|^{\frac{2(n-1)}{n-2}} dx \Bigg)^{\frac{n-2}{n-1}}\\
&\ \ \leq 2 \mathsf{A}(n) \mathcal{E}_{2}(f)
\Bigg( \int_{\R^{n}_{+}} \bigg| \frac{\partial f}{\partial t}(t,x) \bigg|^2 dxdt  \Bigg)^{\frac{1}{2}}\\
&\ \ \leq 2 \mathsf{A}(n) \Bigg(\int_{\R^{n}_{+}} |\nabla_{x}f(t,x)|^2 dxdt\Bigg)^{\frac{1}{2}} \Bigg( \int_{\R^{n}_{+}} \bigg| \frac{\partial f}{\partial t}(t,x) \bigg|^2 dxdt  \Bigg)^{\frac{1}{2}}\\
&\ \ \leq \mathsf{A}(n) \Bigg(\int_{\R^{n}_{+}} |\nabla f(t,x)|^2 dxdt\Bigg).
\end{align*}
So we can say that \eqref{intro3} implies \eqref{intro1}.

For our goal we utilize the weighted $L^p$-affine energy
$$\mathcal{E}_{p}(f, \sigma)=c_{n-1,p} \Bigg( \int_{\sfee}
\|\nabla_{\xi}f\|_{L^p(\R^{n}_+,\sigma)}^{1-n}d\xi\Bigg)^{\frac1{1-n}},$$
where
$$
c_{n,p}=\big(n\omega_{n}\big)^\frac1{n}\Bigg(\frac{n\omega_{n}{\omega_{p-1}}}{2\omega_{n+p-2}}\Bigg)^\frac1p,
$$
$\sigma \colon \R^{n}_+ \to \R_+$ is a measurable function, and $L^p(\R^n_+,\sigma)$ is the weighted Lebesgue space of all measurable functions $f\colon \R^{n}_+ \to \R$ satisfying
$$\|f\|_{L^p(\R^n_+,\sigma)}=\Bigg(\int_{\R^{n}_{+}} |f(t,x)|^p \sigma(t,x)\, dxdt\Bigg)^\frac1p < \infty.$$

We discover an affine non-sharp Poisson trace principle for $\dot{H}^{-\alpha}(\R^{n-1})$ as seen below.

\begin{thm}\label{introT1}
	Let $$
	\begin{cases}
	n\geq 3;\\
	\frac12\le\alpha<1;\\
	p=\frac{2(n-1+2\alpha)}{n+1+2\alpha};\\
	 p'=\frac{p}{p-1}=\frac{2(n+2\alpha-1)}{n+2\alpha-3};\\
	\sigma(t,x)=t^{2\alpha-1}\quad\forall\quad (t,x)\in\R^n_+.
	\end{cases}
	$$ Then for any $g\in C_c^\infty(\R^{n-1})$ and its Poisson extension $$g(t,x)=P_t\ast g(x)\quad\forall\quad (t,x)\in\R^n_+$$ one has
	\begin{equation}\label{intro4}
	\|g\|_{\dot{H}^{-\alpha}(\R^{n-1})}
	<\mathsf{D}(n,p,\alpha) \big(\mathcal{E}_{p}(g,\sigma) \big)^\frac{n-1}{n-1+2\alpha}
	\bigg\|\frac{\partial g}{\partial t}\bigg\|_{L^p(\R^{n}_+, \sigma)}^\frac{2\alpha}{n-1+2\alpha}
	\end{equation}
	where
	 \begin{align*}
	 \mathsf{D}(n,p,\alpha)&=
	 \Bigg(\frac{2^{2\alpha}}{\Gamma(2\alpha)}\Bigg)^{\frac{1}{2}}\left(\pi^{\frac{n-1}{2(1-n-2\alpha)}}(2\alpha)^{\frac{2\alpha}{p(1-n-2\alpha)}}\right)\\ &\ \ \times(n-1)^{-\frac{n-1}{p(n+2\alpha-1)}} 
	 \Bigg( \frac{n+2\alpha-1-p}{p-1}  \Bigg)^{-\frac{1}{q}}\\
	 &\ \ \times\Bigg( \frac{p'\Gamma(\frac{n+1}{2})\Gamma(n+2\alpha-1)}
	 {\Gamma(\frac{2\alpha}{p'})\Gamma(\frac{n-1+p'}{p'})\Gamma(\frac{n+2\alpha-1}{p'})} \Bigg)^{\frac{1}{n+2\alpha-1}}.
	 \end{align*}
	 Consequently, the \eqref{intro4}'s sharp constant
	$$
	\mathsf{D}^\sharp(n,p,\alpha)=\sup_{g\in C^\infty_c(\mathbb R^{n-1})}\left(	\|g\|_{\dot{H}^{-\alpha}(\R^{n-1})} \big(\mathcal{E}_{p}(P_t\ast g,\sigma) \big)^\frac{n-1}{1-n-2\alpha}
	\bigg\|\frac{\partial }{\partial t}P_t\ast g\bigg\|_{L^p(\R^{n}_+, \sigma)}^\frac{2\alpha}{1-n-2\alpha}\right)
	$$
	is strictly less than $\mathsf{D}(n,p,\alpha)$.
\end{thm}

The main ingredient for validating our statement is the coming-up-next sharp affine weighted $L^p$ Sobolev inequality obtained in \cite[Theorem 1.1]{Haddad} by  Haddad, Jim\'enez and Montenegro.

\begin{thm}\label{Haddad THM}
	Let $$
	\begin{cases} a\geq 0;\\
	1\leq p <n+a;\\
	p_{a}^{*}=\frac{p(n+a)}{n-p+a};\\
	\sigma(t,x)=t^a\ \ \forall\ \ (t,x)\in\R^n_+.
	\end{cases}
	$$ 
	Then there is a sharp constant $\mathsf{J}(n,p,a)$ depending only on the triple $(n,p,a)$ such that if $f\in C^\infty_c(\R^n)$ then
	\begin{equation}\label{weighted}
	\|f\|_{L^{p_{a}^{*}}(\R^{n}_{+},\sigma)} \leq \mathsf{J}(n,p,a)\big(\mathcal{E}_{p}(f,\sigma)\big)^{\frac{n-1}{n+a}} 
	\bigg\|\frac{\partial f}{\partial t}\bigg\|_{L^{p}(\R^{n}_{+},\sigma)}^{\frac{1+a}{n+a}}.
	\end{equation}
	The equality in \eqref{weighted} holds if 
	\begin{equation}\label{eq weighted}
	f(t,x)=
	\begin{cases}
	c(1+|\gamma t|^{\frac{p}{p-1}}+ |B(x-x_0)|^{\frac{p}{p-1}})^{\frac{p-n-a}{p}} &\mbox{as}\ p>1;\\
	c\ \chi_{\mathbb{B}^n}(\gamma t, B(x-x_0)) &\mbox{as}\ p=1,
	\end{cases}
	\end{equation}
	for some quadruple
	$$(c,|\gamma|,x_0, B)\in \R\times\R_+\times\R^{n-1}\times GL_{n-1},
	$$ 
	where $\chi_{\mathbb{B}^n}$ is the characteristic function of the unit ball $\mathbb B^n$ in $\R^{n}$.
\end{thm}

We will manipulate the left-hand side of \eqref{weighted} to establish \eqref{intro4}. We observe that also \eqref{weighted} is invariant under the affine transformations of $\R_{+}^n$, and this allows us to say that \eqref{intro4} is affinely invariant. This, along with a non-Poisson viewpoint of \eqref{P1} and the idea of proving Theorem \ref{introT1}, leads us to the following sharp non-Poisson trace principle for $\dot{H}^{-\alpha}(\R^{n-1})$.

 \begin{thm}\label{introT1New}
 	For $$
 	\begin{cases}
 	n\geq 3;\\
 	\frac12\le\alpha<1;\\
 	p=\frac{2(n-1+2\alpha)}{n+1+2\alpha};\\
 	p'=\frac{p}{p-1}=\frac{2(n+2\alpha-1)}{n+2\alpha-3};\\
q=\frac{1+p-n-2\alpha}{p}\\ 
 	\sigma(t,x)=t^{2\alpha-1}\quad\forall\quad (t,x)\in\R^n_+,
 	\end{cases}
 	$$ 
 	let the so-called non-Poisson kernel $$Q_t(x)=Q(t,x)\ \ \forall\ \ (t,x)\in\mathbb R^n_+$$ 
 	be determined by the Fourier transform $\widehat{(\cdot)}=\mathcal{F}{(\cdot)}$ in space:
 	$$
 	\widehat{Q_t}(x)=\frac{\mathcal{F}{\big(1+t^{p'}+|\cdot|^{p'}\big)^q}(x)}{\mathcal{F}{\big(1+|\cdot|^{p'}\big)^q}(x)}\quad\forall\quad (t,x)\in\mathbb R^n_+.
 	$$
 	Then for $h\in C_c^\infty(\R^{n-1})$ equipped with both the anti-Sobole-type norm
 	$$
 	|\|h|\|_{\dot{H}^{-\alpha}(\R^{n-1})}=\left(\int_{\mathbb R^{n-1}}\Big({|x|^{-\alpha}|\widehat{h}(x)|}\Big)^2\bigg(\int_0^\infty\big|s^\alpha\widehat{Q_{\frac{s}{|x|}}}(|x|)\big|^2\frac{ds}{s}\bigg)\,dx\right)^\frac12
 	$$
 	and the non-Poisson type extension $$h(t,x)=Q_t\ast h(x)\quad\forall\quad (t,x)\in\R^n_+$$ one has
 	\begin{align*}
 	\mathsf{D}_\sharp(n,p,\alpha)&=\sup_{h\in C^\infty_c(\mathbb R^{n-1})}\left(	|\|h|\|_{\dot{H}^{-\alpha}(\R^{n-1})} \big(\mathcal{E}_{p}(h,\sigma) \big)^\frac{n-1}{1-n-2\alpha}
 	\bigg\|\frac{\partial h}{\partial t}\bigg\|_{L^p(\R^{n}_+, \sigma)}^\frac{2\alpha}{1-n-2\alpha}\right)\\
 	&=\mathsf{D}(n,p,\alpha) \Big(2^{-2\alpha}{\Gamma(2\alpha)}\Big)^{\frac{1}{2}}\\
 	&=\mathsf{J}(n,p,2\alpha-1).
 	\end{align*}
 \end{thm}

The rest of this article is organized as follows. In \S \ref{s2} we present some essentials about the fractional Sobolev space $\dot{H}^\alpha$ and the fractional Laplacian $(-\Delta)^\frac{\alpha}{2}$, while in \S \ref{s3}\&\S \ref{s4} we verify Theorems \ref{introT1}\&\ref{introT1New}.

\section{Basics on $\dot{H}^\alpha$ and $(-\Delta)^\frac{\alpha}{2}$}\label{s2}

For an integer $n\ge 1$ we define the Fourier transform of a function $f\in L^1(\R^n)$ as 
$$\widehat{f}(x)=\mathcal{F}(f)(x)=\int_{\R^n} f(\xi) \exp\big(-2\pi i \langle x,\xi\rangle\big)\, d\xi.$$
It is well known that we can extend, by Plancherel's theorem, the Fourier transform to every $L^2(\mathbb R^n)$-function.

\begin{defn}
	Given $0< \alpha <1$, the fractional homogeneous Sobolev space $\dot{H}^{\alpha}(\R^n)$ is defined as the completion of all functions $f\in C^{\infty}_{c}(\R^n)$ under the norm $$\|f\|_{\dot{H}^{\alpha}(\R^n)}= \bigg(\int_{\R^n} (2\pi |x|)^{2 \alpha} |\widehat{f}(x)|^2 dx \bigg)^{\frac{1}{2}}.$$
\end{defn}

This definition allows us to link it with the fractional Laplacian $(-\Delta)^{\frac{\alpha}{2}}$. More precisely, if  $\mathcal{F}^{-1}$ is the inverse Fourier transform and
$f\in C^{\infty}_{c}(\R^n)$, then the fractional Laplacian of $f$ is 
$$(-\Delta)^{\frac{\alpha}{2}}f(x)=\mathcal{F}^{-1}\big((2\pi |x|)^{\alpha} \widehat{f}(x)\big),$$
and hence we can say that the fractional Laplacian is an operator enjoying
$$\widehat{(-\Delta)^{\frac{\alpha}{2}}f}(x)=(2\pi |x|)^{\alpha} \widehat{f}(x)$$ and the fractional homogeneous Sobolev norm of $f$ 
$$\|f\|_{\dot{H}^{\alpha}(\R^n)}= \bigg(\int_{\R^n} |(-\Delta)^{\frac{\alpha}{2}}f(x)|^2 dx \bigg)^{\frac{1}{2}},$$ exists as a consequence of Plancherel's theorem.

\begin{proposition}
	$\dot{H}^{\alpha}(\R^n)$ endowed with the inner product 
	$$
	(f,g)_{\dot{H}^{\alpha}(\R^n)}=\int_{\R^n}(2\pi|x|)^\alpha\overline{\widehat{f}(x)}\widehat{g}(x)\,dx
	$$
	is an Hilbert space.
\end{proposition}

We have the following fractional Sobolev inequality; see e.g. \cite{Caffarelli, Cotsiolis, Einav/Loss, Lieb, Xiao} for more details about the fractional Laplacian and its connections with PDE.

\begin{thm}\label{FracSob}
	If 
	$$
	f\in \dot{H}^{0<\alpha<\frac{n}{2}}(\R^n),
	$$
	 then
	\begin{equation}\label{d1.1}
	\Big(\int_{\R^n} |f(x)|^{\frac{2n}{n-2\alpha}}dx \Big)^{\frac{n-2\alpha}{n}}\leq \mathsf{B}(n,\alpha)\Gamma(2(1-\alpha)) 2^{2\alpha-1} \|f\|^{2}_{\dot{H}^{\alpha}(\R^n)},
    \end{equation}
    where $B(n,\alpha)$ is the sharp constant in \eqref{intro2}.
	Moreover, the equality in \eqref{d1.1} holds if and only if $$f(x)=A(\gamma^2 +|x-a|^2)^{\frac{2\alpha-n}{2}}$$
	for some triple 
	$$(A,|\gamma|,a)\in \mathbb{R}\times\R_+\times\R^n.
	$$
\end{thm}

The previous theorem indicates that the fractional homogeneous Sobolev space may be treated as a function space - more precisely - we have the following corollary.
\begin{corollary}
	Under $0< \alpha < \frac{n}{2}$ one has
	$$\dot{H}^{\alpha}(\R^n)=\big\{ f\in L^{\frac{2n}{n-2\alpha}}(\R^n):\ \ \ (-\Delta)^{\frac{\alpha}{2}}f\in L^2(\R^n) \big\} .$$
\end{corollary}

Now, we want to introduce fractional homogeneous Sobolev space and fractional Laplacian with negative exponent: we use the {Riesz potential}. Let $f$ be a smooth function on $\R^n$ with a fast decay at infinity. According to \cite[Section 5.1]{Stein}, under $0<\alpha<\frac{n}{2}$ we define

$$\begin{cases}
	I_{2\alpha}f(x)=\frac{1}{\kappa(2\alpha)}\int_{\R^n}\frac{f(y)}{|x-y|^{n-2\alpha}}dy;\\
	\kappa(2\alpha)=\frac{\pi^{\frac{n}{2}}2^{2\alpha}
		\Gamma(\alpha)}{\Gamma( \frac{n}{2}-\alpha)}.
\end{cases}$$
Since
$I_{2\alpha}f$ is actually the Riesz potential of $f$, we have
$$(-\Delta)^{-\frac{\alpha}{2}}f=I_{\alpha}f.
$$
To justify this last identification, recall the following results (cf. \cite[Section 5.1, Lemma 2]{Stein} and 
\cite[Corollary 5.10]{Lieb/Loss}).

\begin{lemma}\label{Lemma2.8}
	Let $0<\alpha<\frac{n}{2}$.
	\begin{itemize}
		\item[$\rhd$]For every $\phi\in C_{c}^{\infty}(\R^n)$ one has
		$$\int_{\R^{n}} |x|^{2\alpha-n} \phi(x)dx=
		\kappa(2\alpha) \int_{\R^{n}} \frac{\overline{\hat{\phi}(x)}}{(2\pi|x|)^{2\alpha}}dx.$$
		\item[$\rhd$]The Fourier transform of $I_\alpha f$:
		$$\widehat{I_{2\alpha}f}(x)=(2\pi|x|)^{-2\alpha}\widehat{f}(x)$$ holds in the sense of
		$$\int_{\R^{n}} I_{2\alpha}f(x)\overline{g(x)}dx=
		\int_{\R^{n}} \frac{\widehat{f}(x)\overline{\widehat{g}(x)}}{(2\pi|x|)^{2\alpha}}dx\ \  \forall\ \ f,g\in S(\R^n),
		$$ where $S(\R^n)$ is the Schwarz space on $\R^n$.
	\end{itemize}
\end{lemma}

If we take $f=g$ in the last statement of Lemma \ref{Lemma2.8}, then we have
\begin{equation}\label{equation1}
\int_{\R^{n}} (2\pi|x|)^{-2\alpha}|\widehat{f}(x)|^2 dx=
\frac{1}{\kappa(2\alpha)} \int_{\R^{n}} \int_{\R^{n}} \frac{f(x)\overline{f(y)}}{|x-y|^{n-2\alpha}}dxdy.
\end{equation}

We present now the following lemma, a powerful tool to get \eqref{intro2} and \eqref{d1.1}.

\begin{lemma}\label{Lemma HLS}
	Let $\alpha\in (0,\frac{n}{2})$ and $f$, $g\in L^{\frac{2n}{n+2\alpha}}(\R^n)$. Then 
	\begin{align*}
	\Bigg| \int_{\R^{n}} \int_{\R^{n}} \frac{f(x)g(y)}{|x-y|^{n-2\alpha}} dxdy \Bigg|\leq\Bigg( \frac{\pi^{\frac{n}{2}-\alpha}\Gamma(\alpha)}{\Gamma(\frac{n}{2}+\alpha)}\Bigg)\Bigg( \frac{\Gamma(\frac{n}{2})}{\Gamma(n)} \Bigg)^{-\frac{2\alpha}{n}}\|f\|_{L^{\frac{2n}{n+2\alpha}}(\R^n)}\|g\|_{L^{\frac{2n}{n+2\alpha}}(\R^n)}
	\end{align*}
	with equality if and only if 	$$f(x)=c(\gamma^2 + |x-x_0|^2)^{-\frac{n+2\alpha}{2}} $$
	for some triple
	$$
	(c, \gamma,x_0)\in\mathbb{R}\times\R_+\times \R^n
	$$
	and $g$ is equal to $f$ times some constant.
\end{lemma}

Lemma \ref{Lemma HLS} is a particular case of Hardy-Littlewood-Sobolev inequality and we refer to \cite{Lieb/Loss} for its general setting. For our purpose we apply Lemma \ref{Lemma HLS} to \eqref{equation1} and we get exactly
$$\int_{\R^{n}} (2\pi|x|)^{-2\alpha}|\widehat{f}(x)|^2 dx\leq \mathsf{B}(n,\alpha) \Gamma(2(1-\alpha))2^{2\alpha-1}\|f\|^2_{L^{\frac{2n}{n+2\alpha}}(\R^n)},
$$
which is viewed as the dual version of \eqref{intro2} and \eqref{d1.1}. 
In terms of the fractional Laplacian we have
$$\|(-\Delta)^{-\frac{\alpha}{2}}f\|_{L^2(\R^n)}^2\leq\mathsf{B}(n,\alpha) \Gamma(2(1-\alpha))2^{2\alpha-1}\|f\|^2_{L^{\frac{2n}{n+2\alpha}}(\R^n)}, $$
which justifies the introduction of the fractional anti-Sobolev space:
$$\dot{H}^{-\alpha}(\R^n)=\big\{ u\in L^2(\R^n)|\ \exists\ f\in L^{\frac{2n}{n+2\alpha}}(\R^n):\  u=(-\Delta)^{-\frac{\alpha}{2}}f   \big\}.
$$

\section{Proof of Theorem \ref{introT1}}\label{s3}

First of all, we recall Poisson kernel, Poisson extension and their properties that we will use along the proof.

\begin{defn} For $n\in\mathbb N\setminus\{ 1 \}$ and $t\geq 0$ 
	let 
	$P_t \colon \R^{n-1} \to \R$ be the Poisson kernel, i.e., 
	$$P_t(x)={{\Gamma\bigg(\frac{n}{2}\bigg)}\pi^{-\frac{n}{2}}t}{(t^2 + |x|^2)^{-\frac{n}{2}}}\ \ \forall\ \ x\in\R^{n-1}.$$
	Then the Poisson extension or harmonic extension of a function $f\colon \R^{n-1} \to \R$ is the convolution of $P_t$ with $f$:
	$$f(t,x)=P_t \ast f(x)={\Gamma\bigg(\frac{n}{2}\bigg)}{\pi^{-\frac{n}{2}}}
	\int_{\R^{n-1}} {t}{(t^2 + |x-z|^2)^{-\frac{n}{2}}}f(z)\, dz.$$
\end{defn}

\begin{remark} Below is a short list of the Poisson kernel properties (cf. \cite[Chapter 3, Section 2.1]{Stein}).
	\begin{itemize}
		\item[$\rhd$]The reason why $f(t,x)$ is called the harmonic extension of $f$ is that $P_t$ and $f(t,x)$ are harmonics in $\R^{n}_{+}$.
		\item[$\rhd$]
		$$f\in C^{\infty}(\R^{n-1})\Rightarrow f(t,\cdot)\in C^{\infty}(\R^{n-1})\ \forall\ t\in\R_+.
		$$
		\item[$\rhd$]$$\begin{cases} P_t\in L^p(\R^{n-1})\ \ &\forall\ \ \ (t,p)\in\R_+\times [1,\infty);\\ \|P_t\|_{L^1(\R^{n-1})}=1\ \ &\forall\ \ \ t\in\R_+.
		\end{cases}
		$$
		\item[$\rhd$]$$f\in L^p(\R^{n-1})\Rightarrow f(t,\cdot)\in L^p(\R^{n-1})\ \ \forall\ \ t>0.$$	
		\item[$\rhd$]$$\lim_{t \rightarrow 0^+}f(t,x)=f(x)\ \ \ \text{for a.e.}\ \ x\in \R^{n-1}.
		$$
		
		\item[$\rhd$]
		$$\lim_{t \rightarrow 0^+}\big\|f(t,\cdot)-f(\cdot)\big\|_{L^p(\R^{n-1})}=0.$$
		
	\end{itemize}
\end{remark}

Now we employ Theorem \ref{Haddad THM} to verify \eqref{intro4}. As a matter of fact, if
$$
\begin{cases}
a\geq0;\\
1<p<n+a;\\
p_a^*=2;\\
\sigma(t,x)=t^a\quad\forall\quad (t,x)\in\R^n_+,
\end{cases}
$$
 then
$$
p_a^*=\frac{p(n+a)}{n-p+a}=2\ \Leftrightarrow\ p=\frac{2(n+a)}{n+a+2}<n+a.
$$
Let us now consider $g\in C_{c}^{\infty}(\R^{n-1})$, we have 
$$
g\in L^p(\mathbb R^{n-1})\ \forall\ p\geq 1\ \ \&\ \ g(t,x)=P_t \ast g(x).
$$ 
then an application of 
$$P_t\in L^r(\R^{n-1})\quad\forall\quad 1\leq r< \infty
$$ 
and Young's inequality with 
$$r=\frac{n+a}{n+a-1}
$$ 
derives
$$
 g(t,\cdot)\in L^2(\R^{n-1})\quad\forall\quad t\in\R_+.
$$ 
Accordingly, we can evaluate the left-hand side of \eqref{weighted} with $p_a^*=2$, thereby finding (via Fubini's and Plancehrel's theorems)
\begin{align*}
\|g\|^2_{L^2(\R^{n}_{+},\sigma)}&=\int_{0}^{\infty} \int_{\R^{n-1}} |g(t,x)|^2 t^adxdt\\
& =\int_{0}^{\infty} \Big(\int_{\R^{n-1}} |g(t,x)|^2 dx \Big) t^adt\\
&=\int_{0}^{\infty} \Big(\int_{\R^{n-1}} |\widehat{g}(t,x)|^2 dx \Big) t^adt.
\end{align*}
Upon utilizing 
$$g(t,x)=P_t \ast g(x)$$ 
and the convolution property of Fourier's transform, we obtain
 \begin{align*}
 \|g\|^2_{L^2(\R^{n}_{+},\sigma)}&= \int_{0}^{\infty} \Bigg(\int_{\R^{n-1}} |\widehat{P_t}(x)|^2  |\widehat{g}(x)|^2 dx \Bigg) t^a dt\\
 &=\int_{\R^{n-1}} \Bigg(\int_{0}^{\infty} |\widehat{P_t}(x)|^2 t^a dt \Bigg) |\widehat{g}(x)|^2 dx.
 \end{align*}
At the same time, we can take advantage of the explicit expression of the Fourier transform of $P_t$ to evaluate 

$$\int_{0}^{\infty} |\widehat{P_t}(x)|^2 t^a dt$$ via the gamma function - more precisely - a calculation gives
$$\int_{0}^{\infty} |\widehat{P_t}(x)|^2\, t^a dt=\int_0^\infty \Big(\exp\big(-4\pi t|x|\big)\Big)\,t^adt=
\frac{\Gamma(a+1)}{(4\pi|x|)^{a+1}}.$$
This in turn implies
\begin{align*}
\|g\|^2_{L^2(\R^{n}_{+},\sigma)}&=\int_{\R^{n-1}}
\frac{|\widehat{g}(x)|^2 \Gamma(a+1)}{(4\pi|x|)^{a+1}}dx\\
&=\frac{\Gamma(a+1)}{2^{a+1}} \int_{\R^{n-1}}
\frac{|\widehat{g}(x)|^2}{(2\pi|x|)^{a+1}}dx\\
&=
\frac{\Gamma(a+1)}{2^{a+1}} \big\|(-\Delta)^{-\frac{a+1}{4}}g\big\|^2_{L^2(\R^{n-1})}.
\end{align*}
Observe that 
$$
(t,x)\mapsto g(t,x)=P_t\ast g(x)
$$
is in $C^\infty(\mathbb R^n_+)$. So the differential pair $\big\{\frac{\partial g}{\partial t},\nabla_x g\big\}$ exists and the energy pair
$$\left\{\mathcal{E}_p(g,\sigma),\int_{\mathbb R^n}\Bigg|\frac{\partial g(t,x)}{\partial t}\Bigg|^p\,\sigma(t,x)\,dxdt
\right\}
$$
is well-defined. Without loss of generality, we may assume each energy is finite. Then we can approximate $g$ by $C^\infty_c(\mathbb R^n_+)$-functions, and consequently apply Theorem \ref{Haddad THM} to the fractional Laplacian with exponent $-2^{-2}(a+1)$, thereby  obtaining that if 
$$
\begin{cases}
a=2\alpha-1\ge 0;\\
\sigma(t,x)=t^{2\alpha-1};\\
g\in C_c^\infty(\R^{n-1});\\
g(t,x)=P_t\ast g(x)\quad\forall\quad (t,x)\in\R^n_+,
\end{cases}
$$ 
then
\begin{align*}
\|g\|_{\dot{H}^{-\alpha}(\mathbb R^{n-1})}&=\big\|(-\Delta)^{-\frac{a+1}{4}}g\big\|_{L^2(\R^{n-1})}\\
&= 
\Bigg(\frac{2^{a+1}}{\Gamma(a+1)}\Bigg)^{\frac{1}{2}} \|g\|_{L^2(\R^{n}_{+},\sigma)}\\
&\leq \Bigg(\frac{2^{a+1}}{\Gamma(a+1)}\Bigg)^{\frac{1}{2}}\ \mathsf{J}(n,p,a)\big(\mathcal{E}_{p}(g,\sigma)\big)^{\frac{n-1}{n+a
}} 
\bigg\|\frac{\partial g}{\partial t}\bigg\|_{L^{p}(\R^{n}_{+},\sigma)}^{\frac{a+1}{n+a}}.
\end{align*}
Note that $$\mathsf{D}(n,p,\alpha)=\Bigg(\frac{2^{2\alpha}}{\Gamma(2\alpha)}\Bigg)^{\frac{1}{2}}\ \mathsf{J}(n,p,2\alpha-1).
$$ 
So, \eqref{intro4} follows up.

Next, in order to compute $\mathsf{D}(n,p,\alpha)$, recall that
 $\mathsf{J}(n,p,2\alpha-1)$, as the sharp constant in Theorem \ref{Haddad THM} under $\sigma(t,x)=t^{2\alpha-1}$, has been evaluated in \cite{Haddad}. To see this evaluation, for a generic $a\geq 0$ we use \cite[Appendix A]{Haddad} to write $$\mathsf{J}(n,p,a)=\mathsf{L}(n,p,a)\mathsf{M}(n,p,a),$$
where $\mathsf{L}(n,p,a)$ is the best value found in \cite{Cabre, Ros-Oton} concerning a weighted Sobolev inequality for the $L^p$-norm of the gradient - that is -
\begin{align*}
\mathsf{L}(n,p,a)=\Bigg( \frac{(p-1)^{p-1}}{(n+a)(n-p+a)^{p-1}} \Bigg)^{\frac{1}{p}}\Bigg( \frac{2\pi^{\frac{a+2-n}{2}}\Gamma(n+a)\Gamma(\frac{n+a+2}{2})}{\Gamma(\frac{(n+a)(p-1)}{p}+1)\Gamma(\frac{n+a}{p})\Gamma(\frac{1+a}{2})} \Bigg)^{\frac{1}{n+a}},
\end{align*}
while $\mathsf{M}(n,p,a)$ appears in \cite[Lemma 4.1]{Haddad} - namely - if $p'=\frac{p}{p-1}$ then
\begin{align*}
\mathsf{M}(n,p,a)&=(p')^{-\frac{1}{p'}}\pi^{\frac{1-n}{2(n+a)}}(1+a)^{-\frac{1+a}{p(n+a)}}(n-1)^{\frac{1-n}{p(n+a)}}\\
&\ \ \times \Bigg( \frac{n+a}{p} \Bigg)^{\frac{1}{p}}
\Bigg( \frac{p' \Gamma(\frac{n+1}{2}) \Gamma(\frac{n+a+p'}{p'})}{\Gamma(\frac{1+a}{p'})
	\Gamma(\frac{n-1+p'}{p'})}   \Bigg)^{\frac{1}{n+a}}.
\end{align*}
Putting together $\mathsf{L}(n,p,a)$, $\mathsf{M}(n,p,a)$ and $a=2\alpha-1$, one has
\begin{align*}
\mathsf{J}(n,p,2\alpha-1)&=\pi^{-\frac{n-1}{2(n+2\alpha-1)}}(2\alpha)^{-\frac{2\alpha}{p(n+2\alpha-1)}}\\ &\ \ \times(n-1)^{-\frac{n-1}{p(n+2\alpha-1)}} 
\Bigg( \frac{n+2\alpha-1-p}{p-1}  \Bigg)^{-\frac{1}{p'}}\\
&\ \ \times\Bigg( \frac{p'\Gamma(\frac{n+1}{2})\Gamma(n+2\alpha-1)}
{\Gamma(\frac{2\alpha}{p'})\Gamma(\frac{n-1+p'}{p'})\Gamma(\frac{n+2\alpha-1}{p'})} \Bigg)^{\frac{1}{n+2\alpha-1}},
\end{align*}
where
$$
   p=\frac{2(n+2\alpha-1)}{n+2\alpha+1}\quad\&\quad
   p'=\frac{2(n+2\alpha-1)}{n+2\alpha-3}.
$$

Finally, an application of the equality case of Theorem \ref{Haddad THM} derives that the equality in \eqref{intro4} holds if the equality in \eqref{eq weighted} is valid under
$$p=\frac{2(n+a)}{n+a+2}=\frac{2(n+2\alpha-1)}{n+2\alpha+1}.$$
In other words, \eqref{intro4} takes its equality if 
$$g(t,x)=c\Big(1+|\lambda t|^{p'}+
|B(x-x_0)|^{p'}\Big)
^{\frac{p+1-n-2\alpha}{p}}$$
for any given quadruple $$(c,|\lambda|,x_0, B)\in \R\times\R_+\times\R^{n-1}\times GL_{n-1}.
$$ 
But, the last function $g(t,x)$ is never harmonic in $\mathbb R^n$ unless it is a constant. Therefore, the equality of \eqref{intro4} is not valid for any $g\in C_c^\infty(\mathbb R^{n-1})$ - i.e. - \eqref{intro4} is strict. Note that $C^\infty_c(\mathbb R^{n-1})$ is dense in $L^{1\le p<\infty}(\mathbb R^{n-1})$. So, 
$\mathsf{D}^\sharp(n,p,\alpha)$ is smaller than $\mathsf{D}(n,p,\alpha)$.

\section{Proof of Theorem \ref{introT1New}}\label{s4}

On the one hand, we utilize
$$
h(t,x)=Q_t\ast h(x)\ \ \&\ \ \widehat{h(t,\cdot)}(x)=\widehat{Q_t}(x)\widehat{h}(x)
$$
to achieve
\begin{align*}
\|h\|^2_{L^2(\R^{n}_{+},\sigma)}&= \int_{0}^{\infty} \bigg(\int_{\R^{n-1}}|\widehat{Q_t}(x)|^2  |\widehat{h}(x)|^2 dx \bigg) t^{2\alpha -1} dt\\
&=\int_{\R^{n-1}} \bigg(\int_{0}^{\infty}|\widehat{Q_t}(x)|^2 t^{2 \alpha -1}\, dt \bigg) |\widehat{h}(x)|^2 dx\\
&=\int_{\mathbb R^{n-1}}\Big({|x|^{-\alpha}|\widehat{h}(x)|}\Big)^2\bigg(\int_0^\infty\big|s^\alpha\widehat{Q_{\frac{s}{|x|}}}(|x|)\big|^2\frac{ds}{s}\bigg)\,dx\\
&=|\|h|\|^2_{H^{-\alpha}(\mathbb R^{n-1})}.
\end{align*}
Upon using the first part of the argument for Theorem \ref{introT1}, we obtain
\begin{equation}
\label{e41}
\mathsf{D}_\sharp(n,p,\alpha)\le \mathsf{D}(n,p,\alpha) \Bigg(\frac{\Gamma(2\alpha)}{2^{2\alpha}}\Bigg)^{\frac{1}{2}}.
\end{equation}

On the other hand, if 
$$
h_\star(x)=(1+|x|^{p'})^q
$$
then
$$
\begin{cases}
h_\star(t,x)=\big(1+t^{p'}+|x|^{p'}\big)^q=Q_t\ast h_\star(x);\\
h_\star(x)=\lim_{t \rightarrow 0^+}{h_\star(t,x)}.
\end{cases}
$$
In other words, $Q_t(x)$ approaches the Dirac point mass $\delta_x$ at $x$ as $t\to 0^+$ and enjoys the following property
\begin{align*}
\|h_\star(t,\cdot)\|_{L^{1}(\mathbb R^{n-1})}&=\int_{\mathbb R^{n-1}}h_\star(t,x)\,dx\\
&=\int_{\mathbb R^{n-1}}\left(\int_{\mathbb R^{n-1}}Q_t(x-y)\,dx\right)h_\star(y)\,dy\\
&=\|Q_t\|_{L^{1}(\mathbb R^{n-1})}\|h_\star\|_{L^{1}(\mathbb R^{n-1})},
\end{align*}
which in turn implies
$$
\|Q_0\|_{L^{1}(\mathbb R^{n-1})}=\|Q(0,\cdot)\|_{L^{1}(\mathbb R^{n-1})}=\lim_{t \rightarrow 0^+}\|Q_t\|_{L^{1}(\mathbb R^{n-1})}=1.
$$
Moreover, thanks to
$$
h_\star(x)=h_\star(|x|)\ \ \&\ \ h_\star(t,x)= h_\star(t,|x|),
$$
we have
$$
\begin{cases}
\widehat{Q_t}(x)=\widehat{Q_t}(|x|)\ \ \&\ \ Q_t(x)=Q_t(|x|);\\
\int_{0}^{\infty} \big|\widehat{Q_t}(x)\big|^2 t^a dt= \int_{0}^{\infty} \big|\widehat{Q_t}(|x|)\big|^2 t^a dt=|x|^{-2\alpha}\int_{0}^{\infty} \big|\widehat{Q_\frac{s}{|x|}}(|x|)\big|^2 s^{2\alpha-1} ds.
\end{cases}
$$
This actually explains why the anti-Sobole-type norm $|\|h|\|_{\dot{H}^{-\alpha}(\mathbb R^{n-1})}$ is introduced and indeed comparable to the anti-Sobolev norm $\|h\|_{\dot{H}^{-\alpha}(\mathbb R^{n-1})}$ - however - unlike the Poisson function $P$ in \eqref{P1} there is no kernel function $Q$ on $\mathbb R^{n-1}$ such that
$$
Q_t(x)=t^{1-n}Q(t^{-1}x)\quad\forall\quad (t,x)\in\mathbb R^{n}_+.
$$
Accordingly, the equality part of \eqref{weighted} shows that the equality of \eqref{e41} is valid.

\end{document}